\documentclass{amsart}

\usepackage{amssymb}
\usepackage{amsmath}
\usepackage{graphicx}
\usepackage{hyperref}
\usepackage{amsthm}
\usepackage{setspace}
\usepackage{verbatim}
	
\theoremstyle{plain}
\newtheorem{theorem}{Theorem}

\newtheorem{proposition}[theorem]{Proposition}

\theoremstyle{definition}

\newtheorem{example}[theorem]{Example}

\numberwithin{equation}{section}

\newenvironment{Macaulay2}{ \begin{spacing}{0.8} \smallskip}{\smallskip \end{spacing} }

\newcommand{\B}[1]{\mathbb #1}

\newcommand{\F}[1]{\mathfrak #1}

\DeclareMathOperator{\Span}{span}

\newcommand{\LT}{\operatorname{in}_\geq}
\newcommand{\LTg}{\operatorname{in}_\succeq}

\newcommand{\ideal}[1]{\langle #1 \rangle}

\author{Robert Krone}

\title{Numerical Hilbert functions for Macaulay2}

\begin{document}
\begin{abstract}
The {\em NumericalHilbert} package for {\em Macaulay2} includes algorithms for computing local dual spaces of polynomial ideals, and related local combinatorial data about its scheme structure.  These techniques are numerically stable, and can be used with floating point arithmetic over the complex numbers.  They provide a viable alternative in this setting to purely symbolic methods such as standard bases.  In particular, these methods can be used to compute initial ideals, local Hilbert functions and Hilbert regularity.
\end{abstract}

\maketitle

%%%%%%%%%%%%%%%%%%%%%%%%%%%%%%%%%%%%%%%%%%%%%

\section{Introduction}

Symbolic algorithms in algebraic geometry are well studied and can be very powerful, but for many applications they can also be ineffective.  Many problems are most naturally phrased over the real or complex numbers, and the running times of symbolic algorithms may become impractical.  There has been a recent push to develop numerical or symbol-numeric hybrid alternatives.  These methods have already shown success in tackling problems in engineering and other fields \cite{Sommese-Wampler-book-05}, for example using homotopy tracking to quickly solve zero-dimensional polynomial systems \cite{BHSW06}, \cite{Leykin:NAG}.  Numerical methods are good at computing information about varieties, and often do best when things are regular, but we would also like insight into schemes structure.

The Macaulay dual space is a tool that can be used to numerically compute local combinatorial information, such as multiplicity structure and Hilbert function, about a polynomial system at a particular point.  In conjunction with other tools, this information can be used to recover scheme theoretic data about an ideal.  First introduced in \cite{M16}, for a thorough discussion we refer to \cite{MMM96}.  The {\em NumericalHilbert} package for {\em Macaulay2} \cite{M2www} provides methods for both computing dual spaces and extracting information from them, implementing algorithms described in \cite{Mourrain:inverse-systems}, \cite{DZ05} and \cite{Krone:dual-bases-for-pos-dim}.  These methods are designed to work using floating point arithmetic over the complex numbers, and are stable to numerical error.  

Related existing software includes ApaTools by Zeng \cite{zeng2008apatools} and a Maple package of Mantzaflaris and Mourrain \cite{mantzaflaris2011deflation}.  Both can compute multiplicity structure of an isolated solution to a polynomial system, while our package expands this functionality by including tools for computing on positive-dimensional varieties as well.

We will assume we have black-box methods for computing numerical rank, kernel or image of a matrix.  In practice this can be done with singular value decomposition.  Additionally we will assume access to a method that can approximately sample points from a variety, with any desired precision.

\section{Dual space}

Let $R = \B C[x_1,\ldots,x_n]$ and let $R_d$ denote the subspace of polynomials of degree at most $d$.  The {\em truncated dual space} is the vector-space dual $D^d_0 = (R_d)^*$, i.e. the set of functionals $R_d \to \B C$.  The {\em Macaulay dual space} (or {\em local dual space}) of $R$ at the origin is
\[ D_0 = \bigcup_{d \geq 0} D_0^d. \]
For each monomial $x^\alpha \in R$, we define a monomial functional $\partial^\alpha$ by
\[\partial^\alpha\left(\sum_{\beta \in \B N^n} c_\beta x^\beta \right) = c_\alpha. \]
Note that these monomial functionals form a basis of $D_0$, and so its elements can be represented as ``polynomials'' in $\B C[\partial_1,\ldots,\partial_n]$.

For an arbitrary point $y \in \B C^n$ we can also define the dual space at $y$, denoted $D_y$, by translating $y$ to the origin.  The monomial functional $\partial_y^\alpha \in D_y$ evaluates to 1 on $\prod_{i=1}^n (x_i - y_i)^{\alpha_i}$ and to 0 on all other monomials in $x_1-y_1,\ldots,x_n-y_n$.  We will take our point of interest to be the origin without loss of generality, since we can make it so by a change of coordinates.

The classes {\tt PolySpace} and {\tt DualSpace} defined in the {\tt NAGtypes} package are used to represent finite-dimensional subspaces of $R$ and $D_y$ respectively.  {\tt PolySpace} stores a basis of the subspace, and {\tt DualSpace} stores a basis and the base point $y$.  For convenience we express dual space elements as polynomials in $R$ by identifying $\partial^\alpha$ with $x^\alpha$ although this is an abuse of notation.

For an ideal $I \subset R$ there is a subspace of dual functionals $D_0[I] \subset D_0$ which is the orthogonal to $I$.
 \[ D_0[I] = \{p \in D_0 \mid p(f) = 0 \text{ for all }f \in I\}. \]
We refer to $D_0[I]$ as the dual space of $I$, while the truncated dual space of $I$ is $D_0^d[I] = D_0[I] \cap D_0^d$.

The dual space of $I$ exactly characterizes the local properties of $I$ at the origin.  To make this precise, let $R_0$ denote the localization of $R$ at the maximal ideal $\ideal{x_1,\ldots,x_n}$.  Expressing elements of $R_0$ as power series, we can still evaluate them with functionals in $D_0$ and so define $D_0[IR_0]$.
\begin{proposition}
$D_0[I]$ satisfies the following
\begin{itemize}
 \item $D_0[I] = D_0[IR_0].$
 \item The orthogonal to $D_0[I]$ in $R_0$ is exactly $IR_0$.
\end{itemize}
\end{proposition}

Computing a basis for the dual space of $I$ immediately provides useful combinatorial data about the local behavior of $I$ at the origin.  A consequence of the above proposition is that
 \[ \dim_{\B C} D_0[I] = \dim_{\B C} (R_0/IR_0). \]
Therefore if the origin is an isolated solution to the system described by $I$ then its multiplicity is equal to $\dim_{\B C} D_0[I]$.

There is a simple algorithm to compute the dual space of $I$ truncated to degree $d$ from a generating set $F$.  A functional $p \in D_0^d$ is in $D_0^d[I]$ if and only if $p(x^\alpha f) = 0$ for all $f \in F$ and $|\alpha| \leq d$.  Form a matrix with rows the coefficient vectors of each $x^\alpha f$, including only the terms of degree $\leq d$.  The kernel of this matrix consists of the coefficient vectors of $D_0^d[I]$.  For a full discussion of this algorithm we refer to \cite{DZ05}.

The method {\tt truncatedDual} computes a basis of the truncated dual of an ideal $I$ at a point $y$.  If $y$ is known to be an isolated solution then {\tt zeroDimensionalDual} can be used to compute a basis for the full dual space of $I$.  The optional argument {\tt Strategy} allows the user to specify which algorithm to use.  Using {\tt Strategy => DZ} will use the simple algorithm described above.  However the default {\tt Strategy => BM} is an implementation of an algorithm developed by Mourrain in \cite{Mourrain:inverse-systems}, which generally has better asymptotic running time.

\begin{example}
Let $I = \ideal{x_1^2 + x_2^2 - 4, (x_1 - 1)^2} \subset \B C[x_1,x_2]$.  An approximate solution to the system is $y = (1,1.7320508)$.
\begin{Macaulay2}
\begin{verbatim}
i1 : needsPackage "NumericalHilbert";
        
i2 : R = CC[x_1,x_2];

i3 : I = ideal{x_1^2 + x_2^2 - 4, (x_1 - 1)^2}

             2    2       2
o3 = ideal (x  + x  - 4, x  - 2x  + 1)
             1    2       1     1

o3 : Ideal of R

i4 : y = point matrix{{1,1.7320508}};

i5 : zeroDimensionalDual(y, I)

o5 = | -1 .866025x_1-.5x_2 |

o5 : DualSpace

\end{verbatim}
\end{Macaulay2}
Here $D_y[I] = \Span\{-1,\ 0.866025\partial_1-0.5\partial_2\}$ has dimension 2, the multiplicity of the solution approximated by $y$.
\end{example}

\section{Computing Hilbert functions}

If $I$ is positive-dimensional at the point of interest, the dual space of $I$ has infinite dimension and computing a basis is not possible.  However truncated duals can always be computed, the dimensions of which correspond to values of the Hilbert function.
We define the {\em local Hilbert function} of $I$ at 0 to be
 \[ H_{IR_0}(d) := \dim_{\B C} R_0/(IR_0 + \F m^{d+1}) - \dim_{\B C} R_0/(IR_0 + \F m^{d}) \]
where $\F m = \ideal{x_1,\ldots,x_n}$, the maximal ideal at 0.
\begin{proposition}
 $H_{IR_0}(d) = \dim_{\B C}D_0^{d}[I] - \dim_{\B C} D_0^{d-1}[I].$
\end{proposition}
The method {\tt hilbertFunction} can be applied to a {\tt DualSpace} to find its dimension in a particular degree or list of degrees.

\begin{example}\label{ex:HF}
 Consider the cyclic 4-root system $I = \ideal{x_1+x_2+x_3+x_4,\ x_1x_2+x_2x_3+x_3x_4+x_4x_1,\
   x_1x_2x_3+x_2x_3x_4+x_3x_4x_1+x_4x_1x_2,\ x_1x_2x_3x_4-1}$.  It cuts out a curve in $\B C^4$ with several singular points including $(-1,1,1,-1)$.  We obtain an approximation $y$ of this point with our numerical solver.
\begin{Macaulay2}
\begin{verbatim}
i6 : R = CC[x_1..x_4];

i7 : I = ideal {x_1 + x_2 + x_3 + x_4,
                x_1*x_2 + x_2*x_3 + x_3*x_4 + x_4*x_1,
                x_2*x_3*x_4 + x_1*x_3*x_4 + x_1*x_2*x_4 + x_1*x_2*x_3,
                x_1*x_2*x_3*x_4 - 1};

i8 : y = point matrix{{-1.00000000000002-3.43663806114685e-15*ii,
                        .999999999999981+1.03073641806016e-14*ii,
                        1.00000000000002-1.08861573140903e-14*ii,
                       -.999999999999984+3.31351409927520e-15*ii}};

i9 : D = truncatedDual(y, I, 6);

i10 : hilbertFunction({0,1,2,3,4,5,6}, D)

o10 = {1, 2, 1, 1, 1, 1, 1}

o10 : List

\end{verbatim}
\end{Macaulay2}
We see that values of $H_{IR_y}(d)$ for $d$ from 0 to 6 are $1,2,1,1,1,1,1$.  We might guess from this that the Hilbert function remains constant for all $d \geq 2$ but this can't be verified from the information here so far.
\end{example}

For $d >> 0$ the Hilbert function $H_{IR_0}(d)$ is polynomial in $d$ and this is called the {\em Hilbert polynomial}, $HP_{IR_0}(d)$.
If the Hilbert polynomial has degree $r$, then $IR_0$ has dimension $r+1$.  The multiplicity of $I$ at the origin is $cr!$ where $c$ is the lead coefficient of $HP_{IR_0}$.

To compute the Hilbert polynomial of positive-dimensional $IR_0$ we use the dual space to compute the generators of the initial ideal of $IR_0$ using the algorithm described in \cite{Krone:dual-bases-for-pos-dim}.  Let $\succeq$ be a graded monomial order on the monomial functional of $D_0$, and let $\geq$ be the reverse order but on $R_0$, i.e. a graded local order on the monomials.  The initial ideal of $IR_0$, denoted $\LT IR_0$, is the monomial ideal generated by the initial terms of elements of $IR_0$.  Similarly the initial dual space is $\LTg D_0[I]$ spanned by the initial terms of $D_0[I]$.  A monomial $x^\alpha$ is in $\LT IR_0$ if and only if $\partial^\alpha$ is not in $\LTg D_0[I]$, or equivalently
 \[ \LTg D_0[I] = D_0[\LT IR_0]. \]
By observing the monomials missing from $\LTg D_0[I]$ we search for generators of $\LTg D_0[I]$.  The stopping criteria described in \cite{Krone:dual-bases-for-pos-dim} makes this search exhaustive.

The method {\tt gCorners} produces the generators of the initial ideal $\LT IR_0$.  Similarly {\tt sCorners} produces a list of the maximal monomials not in $\LT IR_0$.  To compute the Hilbert polynomial of $IR_0$ we can use the existing method {\tt hilbertPolynomial} on the monomial ideal generated by the output of {\tt gCorners}.  The method {\tt localHilbertRegularity} computes the degree $d$ at which the Hilbert function and Hilbert polynomial agree.

\begin{example}
 We continue with the cyclic 4-root system from Example \ref{ex:HF}.  We omit intermediate output for readability.
\begin{Macaulay2}
\begin{verbatim}
i11 : G = gCorners(y, I)

o11 = | x_4 x_3 x_2^2 x_1x_2 |

              1       4
o11 : Matrix R  <--- R

i12 : hilbertPolynomial(monomialIdeal G, Projective=>false)

o12 = 1

o12 : QQ[i]

\end{verbatim}
\end{Macaulay2}
 The generators of the initial ideal $\LT IR_0$ are $\{x_4, x_3, x_2^2, x_1x_2\}$.  The Hilbert function of an ideal and its initial ideal agree.  The Hilbert polynomial of a monomial ideal such as $\LT IR_0$ is simple to compute and here we see it is $HP_{IR_0}(d) = 1$.  Since the Hilbert polynomial has degree 0, the local dimension of $I$ at $y$ is 1.  The multiplicity is also 1.
\begin{Macaulay2}
\begin{verbatim}
i13 : localHilbertRegularity(y, I)

o13 = 2

\end{verbatim}
\end{Macaulay2}
The regularity of the Hilbert function is $d = 2$.  This confirms that $H_{IR_y}(d) = 1$ for all $d \geq 2$.  We have now completely described the local Hilbert function of $I$ at $y$.
\end{example}

\begin{comment}
\section{Dual space of colon ideals}

Polynomials in $R$ act in a natural way on $D_0$ through a process we refer to as {\em differentiation}.  Define the derivative with respect to $f \in R$ as $\sigma_f:D_0 \to D_0$ which takes a functional $p$ and pre-multiplies it by $f$, so $\sigma_f p(g) = p(fg)$.  For a monomial $x^\alpha$ and a monomial functional $\partial^\beta$, we have $\sigma_{x^\alpha}\partial^\beta = \partial^{\beta - \alpha}$ if $\beta - \alpha$ is non-negative, and 0 otherwise.  To compute general $\sigma_f p$ note that $\sigma$ is bilinear in $f$ and $p$.

Because an ideal $I = \ideal{f_1,\ldots,f_s}$ is closed under multiplication, its dual space $D_0[I]$ is closed under differentiation.  $D_0[I]$ is the maximal subspace of $D_0$ which is closed under differentiation and such that $p(f_i) = 0$ for all $p$ in the subspace and each generator $f_i$ of $I$.  This observation leads to the algorithm in \cite{Mourrain:inverse-systems} for computing $D_0^d[I]$ given a set of generators for $I$, which is the default algorithm used in the package.  For $0 \leq k \leq d$, assume we have a basis for $D_0^{k-1}[I]$.  Then find a basis for the set of $p \in D_0^k$ such that $\sigma_{x_j}p \in D_0^{k-1}[I]$ for $j = 1,\ldots,n$ and $p(f_i) = 0$ for $i = 1,\ldots,s$.

\begin{proposition}
 $\sigma_f D_0[I] = D_0[I:\ideal{f}].$
\end{proposition}
\end{comment}
\bibliographystyle{plain}
\bibliography{bib}

\end{document}